\documentclass{article}
\usepackage{xcolor}
\definecolor{grey}{rgb}{0.5,0.5,0.5}
\usepackage{amsmath}
\usepackage[left=2.6cm, right=2.6cm, bottom=2.6cm, top=2.6cm]{geometry}
\usepackage{hyperref}
\usepackage{cleveref}
\usepackage[dvipsnames]{xcolor}
\usepackage{cancel,xcolor}
\makeatletter

\usepackage{dsfont}
\usepackage{amssymb}
\usepackage{mathtools}
\usepackage{tikz-cd}
\usepackage{amsthm}
\usepackage{graphicx}
\usepackage{graphics}
\usepackage{mathrsfs}
\usepackage{cleveref}
\usepackage{thmtools}
\usepackage{enumitem}
\usepackage{bbm}
\usepackage{dsfont}
\usepackage[toc,page]{appendix}

\newlist{steps}{enumerate}{1}
\usepackage{mathtools}
\makeatletter
\newcommand{\xMapsto}[2][]{\ext@arrow 0599{\Mapstofill@}{#1}{#2}}
\def\Mapstofill@{\arrowfill@{\Mapstochar\Relbar}\Relbar\Rightarrow}
\makeatother
\usepackage{array}
\setlist[steps, 1]{label = Step \arabic*:}
\theoremstyle{plain}
\usepackage{footnote}
\makesavenoteenv{tabular}
\usepackage{amsthm}

\theoremstyle{plain}
\newtheorem{theorem}{Theorem}[section]
\newtheorem{proposition}[theorem]{Proposition}
\newtheorem{lemma}[theorem]{Lemma}
\newtheorem{corollary}[theorem]{Corollary}
\theoremstyle{remark}
\newtheorem{remark}[theorem]{Remark}
\theoremstyle{definition}
\newtheorem{definition}[theorem]{Definition}

\numberwithin{equation}{subsection}
\usepackage{rotating}
\usepackage{comment}
\usepackage[dvipsnames]{xcolor}
\definecolor{darkred}{RGB}{139,0,0}
\definecolor{darkgreen}{RGB}{0,100,0}
\definecolor{darkblue}{RGB}{0,0,139}
\definecolor{Slash1}{RGB}{230,25,75}    
\definecolor{Slash2}{RGB}{60,180,75}    
\definecolor{Slash3}{RGB}{0,130,200}    
\definecolor{Slash4}{RGB}{130,200,200}   
\definecolor{Slash5}{RGB}{0,0,250}   
\definecolor{Slash6}{RGB}{70,240,240}   
\definecolor{Slash7}{RGB}{155,125,25}   
\definecolor{Slash8}{RGB}{210,245,60}   
\definecolor{Slash9}{RGB}{240,180,250}   
\definecolor{Slash10}{RGB}{255,180,195} 
\usepackage{mathabx,epsfig}

\usepackage[nottoc,notlot,notlof]{tocbibind} 

\usepackage[dvipsnames]{xcolor}
\usepackage{cancel}

\definecolor{darkred}{RGB}{40,0,0}
\definecolor{darkgreen}{RGB}{0,80,0}
\definecolor{darkblue}{RGB}{0,0,110}

\title{Isometric Embeddings of Hyperk\"{a}hler $\textup{T}^*\mathbb{CP}^{n}$}
\author{Joshua Lackman}
\date{}
\begin{document}

\maketitle
\begin{abstract}
\noindent 
We construct explicit $U(n)$–equivariant isometric embeddings
\[
\textup{T}^*\mathbb{CP}^{n-1} \xhookrightarrow{}\mathbb{R}^{n^4+2n^2-1}\,,\;\;n\ge 2
\]
with respect to Calabi's hyperk\"{a}hler metric on $\textup{T}^*\mathbb{CP}^{n-1}$. 
To achieve this, we first construct a diffeomorphism from $\textup{T}^*\mathbb{CP}^{n-1}$ onto the manifold of rank–1 projections in $M_n(\mathbb{C}),$ such that the pullback of half of the trace–class norm is a K\"{a}hler potential for the hyperk\"{a}hler metric. We then use the Hilbert space structure of $M_n(\mathbb{C})$ to produce the isometric embeddings into Euclidean space.
\end{abstract}
\tableofcontents
\section{Introduction}
A family of complete hyperk\"{a}hler metrics were obtained by Calabi ([\ref{cal}]) on $\textup{T}^*\mathbb{CP}^{n-1},\,n\ge 2.$ In the special case that $n=2,$ this is the celebrated Eguchi–Hanson space ([\ref{egu}]), which was the first non–flat hyperk\"{a}hler manifold that was discovered. Due to abstract embedding results, it is known that each of these spaces admits an isometric embedding into some Euclidean space, as do all Riemannian manifolds. However, the only explicit isometric embedding of any non–flat hyperk\"{a}hler manifold into Euclidean space that appears in the literature is an embedding of the form\footnote{In [\ref{dun}], the authors construct a local isometric
embedding of the Eguchi–Hanson space into the Minkowski space $\mathbb{R}^{7,1}.$}
\begin{equation}
\textup{T}\mathbb{CP}^1\xhookrightarrow{}\mathbb{R}^{11}\,,
\end{equation}
found in [\ref{han}]. In this paper, we give a simple construction of isometric embeddings
\begin{equation}
    \textup{T}^*\mathbb{CP}^{n-1}\xhookrightarrow{}\mathbb{R}^{n^4+2n^2-1},\;\;n\ge 2
\end{equation}
which are equivariant with respect to the triholomorphic isometry groups, in particular. While Nash–type embedding results give upper bounds on the minimal dimension of a Euclidean space that a Riemannian manifold isometrically embeds into, no such result exists for equivariant Nash–type embeddings ([\ref{mos}]).
\subsection{The Hyperk\"{a}hler Manifold}
Hitchin constructed the hyperk\"{a}hler metric on $\textup{T}^*\mathbb{CP}^{n-1}$ as follows ([\ref{hit}]): let $\mathcal{H}$ be a finite–dimensional complex Hilbert space and let $\mathcal{H},\,\mathcal{H}^*$ be equipped with the natural translation–invariant Riemannian metrics. Let
\begin{equation}
    \mathcal{M}:=\{(v,f)\in\mathcal{H}\times\mathcal{H}^*:f(v)=0,\,\|v\|^2-\|f\|^2=1\}
\end{equation}
be equipped with the pullback of the product metric. There is a $U(1)$–action on $\mathcal{M},$ given by
\begin{equation}
    \lambda\cdot(v,f)=(\lambda v,\bar{\lambda}f)\,,
\end{equation}
and the unique Riemannian metric on $\mathcal{M}/U(1)$ for which $\mathcal{M}\to\mathcal{M}/U(1)$ is a Riemannian submersion is hyperk\"{a}hler. There is a natural identification of $\mathcal{M}/U(1)$ with 
\begin{equation}
 \textup{T}^*\mathbb{P}(\mathcal{H})=\{(v,f)\in \mathcal{H}\oplus\mathcal{H}^*:\|v\|=1,\,f(v)=0\}/U(1)\,,   
\end{equation}
given by 
\begin{equation}
 \textup{T}^*\mathbb{P}(\mathcal{H})\to \mathcal{M}/U(1)\,,\;\;[v,f]\mapsto [\alpha v,\alpha^{-1}f]\,,\;\;\alpha^2=\frac{1+\sqrt{1+4\|f\|^2}}{2}\,,
\end{equation}
and the corresponding hyperk\"{a}hler structure is the one discovered by Calabi.
\subsection{The Isometric Embedding}
Before writing out our isometric embedding, we set up some notation: for a finite–dimensional complex Hilbert space $\mathcal{V}$, $\mathcal{L}(\mathcal{V})$ denotes the Hilbert space of linear maps $\mathcal{V}\to\mathcal{V}$ with the Hilbert–Schmidt inner product, let \begin{equation}
 \vee:\mathcal{V}\to\mathcal{V}^* \,,\;\;x\mapsto x^\vee 
\end{equation}
denote the standard antilinear isomorphism, and for $v\in\mathcal{V},\,f\in\mathcal{V}^*,$ $vf\in\mathcal{L}(\mathcal{V})$ is defined as 
\begin{equation}
  (vf)(u)=f(u)v\,,\;\;u\in\mathcal{V}\,.  
\end{equation}
All Hilbert spaces are equipped with the natural translation–invariant Riemannian metric, ie, the real part of the inner product.
\begin{theorem}\label{main}
The map
\begin{equation}\label{codo}
\mathcal{M}/U(1)\to\mathbb{R}\oplus\mathcal{L}(\mathcal{H})\oplus\mathcal{L}(\mathcal{H}^{\otimes 2})\footnote{We are regarding the right side as a real Hilbert space, while $\mathcal{H}^{\otimes 2}:=\mathcal{H}\otimes_{\mathbb{C}}\mathcal{H}.$}
\end{equation}
given by
\begin{equation}\label{embe}
    [v,f]\mapsto \frac{\sqrt{\|v\|^2+\|f\|^2}}{2}\bigg(\sqrt{\frac{5}{2}},\,\frac{(v+f^\vee)(v^\vee-f)}{\|v\|^2+\|f\|^2},\,\frac{(v+f^\vee)(v^\vee-f)\otimes (v-f^\vee)(v^\vee+f)}{\sqrt{2}(\|v\|^2+\|f\|^2)^2}\bigg)
\end{equation}
is an isometric embedding. Its image is contained in a subspace of dimension $(\textup{dim}_{\mathbb{C}}\mathcal{H})^4+2(\textup{dim}_{\mathbb{C}}\mathcal{H})^2-1.$
\end{theorem}
This isometric embedding is equivariant with respect to $PU(\mathcal{H})\times\mathbb{Z}/2\mathbb{Z},$ with $PU(\mathcal{H})$ being the triholomorphic isometry group of $\mathcal{M}/U(1),$ ([\ref{grg}]). For $A\in PU(\mathcal{H})$ and $[v,f]\in\mathcal{M}/U(1),\,S,T\in\mathcal{L}(\mathcal{H}),$ 
\begin{equation}
    A\cdot [v,f]=[Av,fA^{-1}]\,,\;A\cdot S=ASA^{-1},\;A\cdot S\otimes T=ASA^{-1}\otimes ATA^{-1}\,.
\end{equation}
Meanwhile, the nontrivial element of $\mathbb{Z}/2\mathbb{Z}$ acts by 
\begin{equation}
[v,f]\mapsto [v,-f]\,,\; S\mapsto S^*,\; S\otimes T\mapsto S^*\otimes T^*\,.
\end{equation}
The action on $\mathbb{R}$ is trivial. Furthermore, if $\mathcal{H}=\mathcal{H}_{\mathbb{R}}\otimes_{\mathbb{R}}\mathbb{C}$ for a real Hilbert space $\mathcal{H}_{\mathbb{R}},$ then there is an additional $\mathbb{Z}/2\mathbb{Z}$ isometry given by complex conjugation, and \ref{embe} is still equivariant, making \ref{embe} equivariant with respect to $(PU(\mathcal{H})\rtimes\mathbb{Z}/2\mathbb{Z})\times\mathbb{Z}/2\mathbb{Z},$ where the nontrivial element of $\mathbb{Z}/2\mathbb{Z}$ acts on $PU(\mathcal{H})$ by complex conjugation; with respect to the codomain of \ref{codo}, this additional involution acts trivially on $\mathbb{R}$ and acts by complex conjugation on the other two factors.
\\\\The image of \ref{embe} is contained in the subspace given by
\begin{equation}
    \{(\lambda,A,M):\textup{Im}(\textup{Tr}\,A)=0,\;\sqrt{5}\textup{Tr}\,(\Sigma M)=\lambda,\;M=\Sigma M^*\Sigma\}\,,
\end{equation}
where $\Sigma\in\mathcal{L}(\mathcal{H}^{\otimes 2})$ is the canonical flip, ie, $v\otimes w\mapsto w\otimes v.$
\subsection{The Construction of the Embedding}
A more illuminating way of expressing this embedding is as a composition of three maps: the first map is a diffeomorphism
\begin{equation}
    \mathcal{M}/U(1)\to \textup{rank–1 idempotents}\in\mathcal{L}(\mathcal{H})\,,\;\;[v,f]\mapsto (v+f^\vee)(v^\vee-f)\,;
\end{equation}
the second map is of the form 
\begin{equation}
    \mathcal{V}-\{0\}\xhookrightarrow{} \mathbb{R}\oplus \mathcal{V}\oplus\mathcal{L}(\mathcal{V})_{\textup{self–adjoint}}
\end{equation}
for an arbitrary finite–dimensional Hilbert space $\mathcal{V},$ where in \ref{embe}, $\mathcal{V}=\mathcal{L}(\mathcal{H});$ the third map uses the unitary map given by
\begin{equation}\label{uni}
   \mathcal{L}(\mathcal{H})^*\to \mathcal{L}(\mathcal{H})\,,\;\;f\mapsto (f^{\vee})^*\,.
\end{equation}
Noting that $\mathcal{L}(\mathcal{L}(\mathcal{H}))\cong \mathcal{L}(\mathcal{H})\otimes \mathcal{L}(\mathcal{H})^*,$ from \ref{uni} we obtain a unitary map $\mathcal{L}(\mathcal{L}(\mathcal{H}))\to \mathcal{L}(\mathcal{H}^{\otimes 2}).$
\section{Hyperk\"{a}hler Structure of $\textup{T}^*\mathbb{CP}^{n-1}$}
First, we will present the standard hyperk\"{a}hler quotient construction of $\textup{T}^*\mathbb{CP}^{n-1},$ ([\ref{hit}]). We will then give another perspective on it using projections.
\subsection{Hyperk\"{a}hler Quotient Construction}
Let $(\mathcal{H},\langle\cdot,\cdot\rangle)$ be a finite–dimensional complex Hilbert space and let $\vee:\mathcal{H}\to\mathcal{H}^*$ denote the standard antilinear isomorphism. We use the same notation to denote its inverse. The space 
\begin{equation}
\mathcal{H}\oplus\mathcal{H}^*\,.
\end{equation}
is a hyperk\"{a}hler vector space with respect to the product metric
\begin{equation}\label{met}
    \langle(v,f),(w,g)\rangle=\textup{Re}(\langle v,w\rangle+\langle f,g\rangle)
\end{equation}
and almost complex structures
\begin{align}
    &I(v,f)=(iv,if)\,,
    \\& J(v,f)=(f^\vee,-v^\vee)\,,
    \\& K(v,f)=(if^\vee,-iv^\vee)\,.
\end{align}
\begin{remark}
With respect to the hyperk\"{a}hler quotient construction of $\textup{T}^*\mathbb{CP}^{n-1},$ $I$ descends to the natural almost complex structure on the cotangent bundle, whereas $K$ descends to an almost complex structure that we will identify with the affine almost complex structure on the space of rank–1 idempotents.
\end{remark}
\begin{lemma}([\ref{kos}])
A K\"{a}hler potential, for all three almost complex structures (ie, a hyperk\"{a}hler potential), is given by 
\begin{equation}\label{kpot}
  \mathcal{H}\oplus\mathcal{H}^*\to\mathbb{R}\,,\;\;  (v,f)\mapsto \frac{1}{2}(\|v\|^2+\|f\|^2)\,.
\end{equation}   
\end{lemma}
That is, applying any of $\frac{1}{2}IdId,\,\frac{1}{2}JdJd,\,\frac{1}{2}KdKd$ to \ref{kpot} results in \ref{met}, where $d$ is the exterior derivative and $I,J,K$ act in the left–most argument.
\\\\The unitary group $U(1)$ acts on $\mathcal{H}\oplus\mathcal{H}^*$ by $\lambda\cdot(v,f)=(\lambda v,\bar{\lambda} f).$ Under the standard isomorphism of $\frak{u}(1)^*\cong\mathbb{R},$ the hyperk\"{a}hler moment map is given by
\begin{equation}
  \mu:  \mathcal{H}\oplus\mathcal{H}^*\to \mathbb{R}\oplus\mathbb{C}\,,\;\mu(v,f)=\big((\|v\|^2-\|f\|^2-1)/2,f(v)\big)\,.
\end{equation} 
The preimage of $0$ is given by
\begin{equation}
    \{(v,f):\|v\|^2-\|f\|^2=1,\;f(v)=0\}\,,
\end{equation}
and therefore, the hyperk\"{a}hler quotient is given by
\begin{equation}\label{quo}
    (\mathcal{H}\oplus\mathcal{H}^*)/\!/U(1)=\{(v,f)\in \mathcal{H}\oplus\mathcal{H}^*:\|v\|^2-\|f\|^2=1,\;f(v)=0\}/U(1)\,.
\end{equation}
\begin{proposition}\label{pullk}
The pullback of the K\"{a}hler potential \ref{kpot} to $\mu^{-1}(0)$ descends to a K\"{a}hler potential on \ref{quo}, with respect to $J,K.$ 
\end{proposition}
\begin{proof}
Let $\phi:\mu^{-1}(0)\to\mathbb{R}$ be defined by $\phi(v,f)=(\|v\|^2+\|f\|^2)/2.$ We just need to show that for any vertical vector $\xi,$ $d\phi(J\xi)=d\phi(K\xi)=0.$ We have that 
\begin{equation}
    d\phi_{(v,f)}=\textup{Re}(\langle v,\cdot\rangle+\langle f,\cdot\rangle)\,.
\end{equation}
Now, consider an arbitrary vertical vector $\xi_{(v,f)}=(ixv,-ixf),\,x\in \mathbb{R}.$ It follows that 
\begin{equation}
J(\xi_{(v,f)})=(ixf^\vee,ixv^\vee)\,,\;\;K(\xi_{(v,f)})=(-xf^\vee,-xv^\vee)\,.    
\end{equation}
Since $f(v)=0,$ the result follows.
\end{proof}
\subsection{The Corresponding Hyperk\"{a}hler Structure on the Space of Projections}
We will now give an alternative perspective on the aforementioned hyperk\"{a}hler structure. Let $ (v,f)\in\mathcal{H}\oplus\mathcal{H}^*$ satisfy
\begin{equation}
  \|v\|^2-\|f\|^2=1,\;f(v)=0\,. 
\end{equation}
Then
\begin{equation}
    (v+f^\vee)(v^\vee-f)\in\mathcal{L}(\mathcal{H})
\end{equation}
is an idempotent, and by cyclicity of the trace, it has rank–$1.$ This motivates the following:
\begin{definition}
We define a map
\begin{equation}\label{ident}
\Phi:(\mathcal{H}\oplus\mathcal{H}^*)/\!/U(1)\to  \mathcal{L}(\mathcal{H}) \,,\;\;   \Phi([v,f])=(v+f^\vee)(v^\vee-f)\,.
\end{equation}
\end{definition}
\begin{proposition}
$\Phi$ is a holomorphic embedding with respect to $K$ and the natural complex structure on $\mathcal{L}(\mathcal{H}).$ Furthermore, its image is the space of rank–1 idempotents and it intertwines $[v,f]\mapsto [v,-f]$ with the adjoint map.  
\end{proposition} 
\begin{proof}
To see that $\Phi$ is holomorphic, it is enough to check that the map
\begin{equation}\label{varphi}
   \varphi: \mathcal{H}\oplus\mathcal{H}^*\to \mathcal{L}(\mathcal{H})\,,\;\;(v,f)\mapsto (v+f^\vee)(v^\vee-f)
\end{equation}
is holomorphic. Letting $(X,Y)\in\mathcal{H}\oplus \mathcal{H}^*,$ we have that 
\begin{equation}\label{tr}
    d\varphi_{(v,f)}(X,Y)=(X+Y^\vee)(v^\vee-f)+(v+f^\vee)(X^\vee-Y)\,.
\end{equation}
From this, we get that 
\begin{align}
    d\varphi_{(v,f)}(K(X,Y))&=d\varphi_{(v,f)}(iY^\vee,-iX^\vee)\\&=(iY^\vee+iX)(v^\vee-f)+(v+f^\vee)(-iY+iX^\vee)
    =id\varphi_{(v,f)}(X,Y) \,.
\end{align}
To see that $\Phi$ is surjective, we note that if $q$ is a rank–1 idempotent, then $q=uw^\vee$ for some $u,w\in \mathcal{H}$ such that $u$ generates the image of $q,$ $w$ generates the image of $q^*,$ $\langle w,u\rangle =1$ and $\|u\|=\|w\|.$ We then let 
\begin{equation}
 v=(u+w)/2\,,\;\;f=(u^\vee-w^\vee)/2\,.   
\end{equation}
Then $\|v\|^2-\|f\|^2=1,\,f(v)=0.$ To see injectivity, we note that $q$ determines such a $u,w$ up to the $U(1)$–action, and therefore, determines such an $v,f$ up to the $U(1)$–action. 
\\\\To see that $\Phi$ is an immersion, suppose that $d\varphi_{(v,f)}(X,Y)=0.$ From \cref{varphi}, it follows that 
\begin{equation}
    X+Y^\vee=(x+iy)(v+f^\vee)\,,\;\;X^\vee-Y=-(x+iy)(v^\vee-f)
\end{equation}
for some $x,y\in\mathbb{R}.$ Solving for $X,Y,$ we get that 
\begin{equation}\label{solve}
    X=xf^\vee+iyv\,,\;\;Y=xv^\vee-iyf\,.
\end{equation}
On the other hand, 
\begin{equation}
 (X,Y)\in \textup{T}_{(v,f)}\{(w,g)\in\mathcal{H}\oplus\mathcal{H}^*:g(w)=0\}  
\end{equation}
if and only if $Y(v)+f(X)=0.$ Assuming this equation holds and combining it with \cref{solve}, we get that 
\begin{equation}
    x(\|v\|^2+\|f\|^2)=0\,.
\end{equation}
If $\|v\|^2-\|f\|^2=1,$ then it must be the case that $x=0.$ Therefore, 
\begin{equation}
    (X,Y)=(iyv,-iyf)\,,
\end{equation}
ie, $(X,Y)$ is tangent to the $U(1)$–orbit. This completes the proof.
\end{proof}
\begin{lemma}
Let $(v,f)\in (\mathcal{H}\oplus\mathcal{H}^*)/\!/U(1).$ Then $\|v\|^2+\|f\|^2$ equals the Hilbert–Schmidt norm of $\Phi([v,f]).$
\end{lemma}
\begin{proof}
\begin{align}
\Phi([v,f])^*\Phi([v,f])&=(v-f^\vee)(v^\vee+f)(v+f^\vee)(v^\vee-f)
\\&=(v-f^\vee)(\|v\|^2+\|f\|^2)(v^\vee-f)=(\|v\|^2+\|f\|^2)(v-f^\vee)(v^\vee-f)\,.
\end{align}
Taking the trace, we get $(\|v\|^2+\|f\|^2)^2,$ and this completes the proof.
\end{proof}
\begin{remark}
For rank–1 operators, the Hilbert–Schmidt norm equals the trace–class norm. If one constructs the hyperk\"{a}hler structure of the cotangent bundle of a Grassmanian using the analogous hyperk\"{a}hler quotient construction, then the K\"{a}hler potential is the trace–class norm divided by $2.$
\end{remark}
As we will see in the next section, for any finite–dimensional complex Hilbert space $\mathcal{V},$ the norm is a K\"{a}hler potential for a Riemannian metric on $\mathcal{V}-\{0\}.$
\begin{corollary}\label{isoc}
$\Phi$ is an isometric embedding into $\mathcal{L}(\mathcal{H})-\{0\},$ where the latter is equipped with the Riemannian metric whose K\"{a}hler potential is given by 
\begin{equation}
\mathcal{L}(\mathcal{H})-\{0\}\to\mathbb{R}\,,\;\; A\mapsto \frac{\sqrt{\textup{Tr}(A^*A)}}{2}\,.
\end{equation}
\end{corollary}
\section{Isometric Embeddings}
\subsection{Hilbert Space Embedding}
 Due to \cref{isoc}, we have an isometric embedding
\begin{equation}
    \textup{T}^*\mathbb{CP}^{n-1}\xhookrightarrow{} \textup{M}_n(\mathbb{C})-\{0\}\,,
\end{equation} 
where the latter is equipped with the Riemannian metric whose K\"{a}hler potential is half of the Hilbert–Schmidt norm. Therefore, we can obtain an isometric embedding of $\textup{T}^*\mathbb{CP}^{n-1}$ into Euclidean space by finding an isometric embedding of $\mathcal{V}-\{0\}$ into Euclidean space, where $\mathcal{V}$ is an arbitrary Hilbert space and $\mathcal{V}-\{0\}$ is equipped with the Riemannian metric whose K\"{a}hler potential is half of the norm.
\begin{proposition}\label{kap}
Let $\mathcal{V}$ be a finite–dimensional complex Hilbert space. Then 
\begin{equation}
    \mathcal{V}-\{0\}\to\mathbb{R}\,,\;\;x\mapsto\frac{\|x\|}{2}
\end{equation}
is a K\"{a}hler potential\footnote{The convention we are using is that $\tau$ is a K\"{a}hler potential for $g$ if $g=\frac{1}{2}IdId\tau.$} for the Riemannian metric $g$ given by
\begin{equation}\label{h0}
 g_x(v,w):=\textup{Re}\Big(\frac{\langle v,w\rangle}{2\|x\|}-\frac{1}{4\|x\|^3}\langle v,x\rangle\langle x,w\rangle\Big)\,, \;\;x\in\mathcal{V}-\{0\},\,v,w\in\mathcal{V}\,.\footnote{In this definition, we are using the natural trivialization of $\textup{T}(\mathcal{V}-\{0\}).$}
\end{equation}
\end{proposition}
The Cauchy–Schwarz inequality implies that $g$ is positive definite: 
\begin{equation}
    \|v\|^2\|x\|^2\ge |\langle v,x\rangle|^2\implies \frac{\|v\|^2}{\|x\|}\ge \frac{|\langle v,x\rangle|^2}{\|x\|^3}\,.
\end{equation}
Hence, $g$ is, indeed, a Riemannian metric.
\begin{proof}(of \cref{kap})
Let $r:\mathcal{V}-\{0\}\to\mathbb{R}$ be the norm. We will compute $dd^cr(\cdot,i\cdot),$ where
$d^c r:=-dr\circ i.$
We have
\begin{equation}
dr_x(v)=\frac{\operatorname{Re}\langle x,v\rangle}{\|x\|}\,,
\end{equation}
and differentiating again gives
\begin{equation}
\operatorname{Hess}_x r(v,w)
=
\frac{\operatorname{Re}\langle v,w\rangle}{\|x\|}
-
\frac{
\operatorname{Re}\langle x,v\rangle
\operatorname{Re}\langle x,w\rangle
}{\|x\|^3}.
\end{equation}
It follows that 
\begin{equation}
dd^c r(v,Iw)
=
\operatorname{Hess}r(v,w)
+
\operatorname{Hess}r(Iv,Iw)\,,
\end{equation}
and therefore,
\begin{equation}
\begin{aligned}
dd^c r(v,Iw)
=
\frac{2\operatorname{Re}\langle v,w\rangle}{\|x\|}
-\frac{
\operatorname{Re}\langle x,v\rangle
\operatorname{Re}\langle x,w\rangle
+
\operatorname{Re}\langle x,Iv\rangle
\operatorname{Re}\langle x,Iw\rangle
}{\|x\|^3}\,.
\end{aligned}
\end{equation}
Since
\begin{equation}
\operatorname{Re}\langle x,v\rangle
\operatorname{Re}\langle x,w\rangle
+
\operatorname{Re}\langle x,Iv\rangle
\operatorname{Re}\langle x,Iw\rangle
=
\operatorname{Re}
\bigl(
\langle v,x\rangle\langle x,w\rangle
\bigr)\,,
\end{equation}
we obtain
\begin{equation}
dd^c r(v,Iw)
=
\operatorname{Re}\left(
\frac{2}{\|x\|}\langle v,w\rangle
-
\frac{1}{\|x\|^3}
\langle v,x\rangle\langle x,w\rangle
\right)\,.
\end{equation}
Dividing by four completes the proof.
\end{proof}
\begin{lemma}\label{hemb}
Let $\mathcal{V}$ be a finite–dimensional complex Hilbert space and let $\mathcal{V}-\{0\}$ be equipped with the Riemannian metric given by the real part of \ref{h0}. Let
\begin{equation}
\mathbb{R}\oplus\mathcal{V}\oplus\mathcal{L}(\mathcal{V})
\end{equation}
be equipped with the products of the translation invariant Riemannian metrics.
\\\\Then  
\begin{equation}\label{hembb}
\mathcal{V}-\{0\}\to \mathbb{R}\oplus\mathcal{V}\oplus\mathcal{L}(\mathcal{V})\,,\;\;x\mapsto \frac{\sqrt{\|x\|}}{2}\bigg(\sqrt{\frac{5}{2}},\|x\|^{-1}x,\frac{\|x\|^{-2}}{\sqrt{2}}xx^\vee\bigg)
\end{equation}
is an isometric embedding, where for $y\in\mathcal{V},$ $xx^\vee(y):=\langle x,y\rangle x.$
\end{lemma}
\begin{proof}
In the following, $x\in\mathcal{V}-\{0\},\,v,w\in\mathcal{V}.$ 
\\\\The derivatives of $x\mapsto \|x\|^{1/2},\,x\mapsto \|x\|^{-1/2}x,\,x\mapsto \|x\|^{-3/2}xx^\vee$ are, respectively,
\begin{align}
    &v\mapsto \frac{1}{2}\|x\|^{-3/2}\textup{Re}\langle x,v\rangle\,,
    \\& v\mapsto \|x\|^{-1/2}v-\frac{1}{2}\|x\|^{-5/2}\textup{Re}\langle x,v\rangle x\,,
    \\&v\mapsto -\frac{3}{2}\|x\|^{-7/2}\textup{Re}\langle x,v\rangle xx^\vee+\|x\|^{-3/2}(vx^\vee+xv^\vee)\,.
\end{align}
It follows that the pullback metric is one quarter of 
\begin{equation}
(v,w)\mapsto 2\|x\|^{-1}\textup{Re}\langle v,w\rangle-2\|x\|^{-3}\textup{Re}\langle x,v\rangle \textup{Re}\langle x,w\rangle  +\|x\|^{-3}\textup{Re}(\langle x,v\rangle\langle x,w\rangle)\,.
\end{equation}
Since $\textup{Re}(\overline{a}b)+\textup{Re}(ab)=2\textup{Re}(a)\textup{Re}(b)$ for $a,b\in\mathbb{C},$ the result is one quarter of
\begin{equation}
(v,w)\mapsto 2\|x\|^{-1}\textup{Re}\langle v,w\rangle-\|x\|^{-3}\textup{Re}(\langle v,x\rangle\langle x,w\rangle)
\end{equation}
and this completes the proof.
\end{proof}
The image of \ref{hembb} is contained in
\begin{equation}
    \{(\lambda,v,A)\in \mathbb{R}\oplus\mathcal{V}\oplus\mathcal{L}(\mathcal{V}):A=A^*,\,\sqrt{5}\,\textup{Tr}(A)=\lambda\}\,,
\end{equation}
which is a \textit{real} linear subspace of dimension $(\textup{dim}_{\mathbb{C}}\,\mathcal{V})^2+2\,\textup{dim}_{\mathbb{C}}\,\mathcal{V}.$ As a corollary, letting $\mathcal{V}:=\mathcal{L}(\mathcal{H}),$ we get an isometric embedding
\begin{equation}
    \textup{T}^*\mathbb{CP}^{n-1}\xhookrightarrow{}\mathbb{R}^{n^4+2n^2-1}.
\end{equation}
We still consider the case that $\mathcal{V}:=\mathcal{L}(\mathcal{H}).$ As mentioned earlier, we can use the natural unitary map 
\begin{equation}
    \mathcal{L}(\mathcal{H})^*\to\mathcal{L}(\mathcal{H})\,,\;\;f\mapsto (f^\vee)^*
\end{equation}
to replace $\mathcal{L}(\mathcal{L}(\mathcal{H}))$ in \ref{hembb} with $\mathcal{L}(\mathcal{H}^{\otimes 2}),$ and the result is \cref{main}. In particular, we can equivalently write the embedding \ref{hembb} as 
\begin{equation}
    \mathcal{L}(\mathcal{H})-\{0\}\xhookrightarrow{} \mathbb{R}\oplus \mathcal{L}(\mathcal{H})\oplus \mathcal{L}(\mathcal{H}^{\otimes 2})\,,\;\;A\mapsto \frac{\sqrt{\|A\|}}{2}\bigg(\sqrt{\frac{5}{2}},\frac{A}{\|A\|},\frac{1}{\sqrt{2}\|A\|^2}A\otimes A^*\bigg)\,.
\end{equation}
\section*{Acknowledgments}
Thank you to Omar Kidwai for your superb accommodations. A special thanks to Dr. Zhang for all of your support.  
\section{References}
\begin{enumerate}
   \item\label{cal}
Eugenio Calabi. M\'etriques k\"ahl\'eriennes et fibr\'es holomorphes.
Annales scientifiques de l'\'Ecole Normale Sup\'erieure 12 (1979): 269--294.
DOI: 10.24033/asens.1367.

\item\label{dun}
Maciej Dunajski and Paul Tod. Conformal and isometric embeddings of gravitational instantons.
Geometry, Lie Theory and Applications 16 (2022): 21--48.
DOI: 10.1007/978-3-030-81296-6\_2.

\item\label{egu}
Tohru Eguchi and Andrew J. Hanson. Asymptotically flat self-dual solutions to Euclidean gravity.
Physics Letters B 74 (1978): 249--251.
DOI: 10.1016/0370-2693(78)90566-X.

\item\label{grg}
G. W. Gibbons, P. Rychenkova, and R. Goto.
HyperK\"ahler quotient construction of BPS monopole moduli spaces.
Communications in Mathematical Physics 186 (1997): 581--599.
DOI: 10.1007/s002200050121.

\item\label{han}
Andrew J. Hanson and Ji-Ping Sha. Isometric Embedding of the A1 Gravitational Instanton.
Memorial Volume for Kerson Huang (2017): 95--111.
DOI: 10.1142/9789813207431\_0013.

\item\label{hit}
Nigel J. Hitchin. Metrics on moduli spaces.
The Lefschetz Centennial Conference, Part I,
Contemporary Mathematics 58, Part I (1986): 157--178.
DOI: 10.1090/conm/058.1/860411.

\item\label{kos}
Piotr Kobak and Andrew Swann.
HyperK\"ahler potentials via finite-dimensional quotients.
Geometriae Dedicata 88 (2001): 1--19.
DOI: 10.1023/A:1013174027792.

\item\label{mos}
John Douglas Moore and Roger Schlafly.
On equivariant isometric embeddings.
Mathematische Zeitschrift 173 (1980): 119--133.
DOI: 10.1007/BF01159954.

\end{enumerate}
\bigskip
\noindent
Department of Mathematics, Beijing Institute of Technology, Beijing, China\\
E-mail: \texttt{joshualackman@outlook.com}
\end{document}